\documentclass[preprint,1p]{elsarticle}

\makeatletter
\def\ps@pprintTitle{%
	\let\@oddhead\@empty
	\let\@evenhead\@empty
	\def\@oddfoot{\footnotesize\itshape
		{} \hfill}%
	\let\@evenfoot\@oddfoot
}
\makeatother
\usepackage[unicode]{hyperref}

\usepackage{latexsym}
\usepackage{indentfirst}
\usepackage{amsxtra}
\usepackage{amssymb}
\usepackage{amsthm}
\usepackage{amsmath}
\usepackage{mathrsfs} 

\usepackage[capitalise]{cleveref}

\newtheorem{theor}{Theorem}
\newtheorem{prop}[theor]{Proposition}

\newtheorem{cor}[theor]{Corollary}
\theoremstyle{definition}               
\newtheorem{defin}[theor]{Definition}
\newtheorem{ex}{Example}
\newtheorem{exs}[ex]{Examples}
\newtheorem{rem}[theor]{Remark}

\DeclareMathOperator{\id}{id}

\DeclareMathOperator{\im}{Im}
\DeclareMathOperator{\ind}{i}
\DeclareMathOperator{\per}{p}

\newcommand{\alphaa}[2]{\alpha_{#1}{#2}}
\newcommand{\betaa}[2]{\beta_{#1}{#2}}
\newcommand{\thetaa}[2]{\theta_{#1}{#2}}
\newcommand{\indd}[1]{\ind{#1}}
\newcommand{\perr}[1]{\per{#1}}

\begin{document}
\begin{frontmatter}
	\title{Set-theoretical solutions of the Yang-Baxter and pentagon equations on semigroups\tnoteref{mytitlenote}}
	\tnotetext[mytitlenote]{{This work was partially supported by the Dipartimento di Matematica e Fisica ``Ennio De Giorgi" - Universit\`{a} del
			Salento. The authors are members of GNSAGA (INdAM).}}
	\author [unile] {Francesco~CATINO\corref{cor1}}
	\ead{francesco.catino@unisalento.it}
	\cortext[cor1]{Corresponding author}
	\author [unile] {Marzia~MAZZOTTA}
	\ead{marzia.mazzotta@unisalento.it}
	\author [unile] {Paola~STEFANELLI}
	\ead{paola.stefanelli@unisalento.it}
	\address[unile]{Dipartimento di Matematica e Fisica ``Ennio De Giorgi"
		\\
		Universit\`{a} del Salento\\
		Via Provinciale Lecce-Arnesano \\
		73100 Lecce (Italy)\\}
	\begin{abstract}
	\small{The Yang-Baxter and pentagon equations are two well-known equations of Mathematical Physic.
	If $S$ is a set, a map $s:S\times S\to S\times S$ is said to be a \emph{set theoretical solution of the Yang-Baxter equation} if
	$$
	    s_{23}\, s_{13}\, s_{12}
	    = s_{12}\, s_{13}\, s_{23},
	 $$
	 where $s_{12}=s\times \id_S$,  $s_{23}=\id_S\times s$, and $s_{13}=(\id_S\times \tau)\,s_{12}\,(\id_S\times \tau)$ and $\tau$ is the flip map, i.e., the map on $S\times S$ given by $\tau(x,y)=(y,x)$. Instead, $s$ is called a \emph{set-theoretical solution of the pentagon equation} if
	 $$
	    s_{23}\, s_{13}\, s_{12}=s_{12}\, s_{23}.
	 $$
	 The main aim of this work is to display how solutions of the pentagon equation turn out to be a useful tool to obtain new solutions of the Yang-Baxter equation.
	 Specifically, we present a new construction of solutions of the Yang-Baxter equation involving two specific solutions of the pentagon equation. 
	 To this end, we provide a method to obtain solutions of the pentagon equation on the matched product of two semigroups, that is a semigroup including the classical Zappa product.}
	\end{abstract} 
	

\end{frontmatter}

\section{Introduction}
Let $S$ be a set. A \emph{set-theoretical solution of the pentagon equation} on $S$ is a map $s:S\times S\to S\times S$ such that
\begin{equation*}
s_{23}\, s_{13}\, s_{12}=s_{12}\, s_{23},
\end{equation*}
where $s_{12}=s\times \id_S$,  $s_{23}=\id_S\times s$, and $s_{13}=(\id_S\times \tau)\,s_{12}\,(\id_S\times \tau)$ and $\tau$ is the flip map, i.e., the map on $S\times S$ given by $\tau(x,y)=(y,x)$.\\
Set-theoretical solutions of the pentagon equation are used in the pioneering work of Baaj and Skandalis \cite{BaSk93} to obtain multiplicative unitary operators on a Hilbert space.
Later, there have been appeared several papers on this topic: Zakrzewski \cite{Za92}, Baaj and Skandalis \cite{BaSk98}, Kashaev and Sergeev \cite{KaSe98}, Jiang and Liu \cite{JiLi05}, Kashaev and Reshetikhin \cite{KaRe07}, Kashaev \cite{Ka11}, and Catino, Mazzotta, and Miccoli \cite{CaMaMi19x}. \\
On the other hand, a map $s:S\times S\to S\times S$ is a \emph{set-theoretical solution of the quantum Yang-Baxter equation} \cite{Dr92} on a set $S$ if the relation
\begin{equation*}
s_{23}\, s_{13}\, s_{12}
=
s_{12}\, s_{13}\, s_{23}
\end{equation*}
is satisfied, with the same notation adopted for the pentagon equation. Just to give an example, if $f, g$ are idempotent maps from $S$ into itself such that $fg=gf$, then the map $s$ defined as $s(x,y)=(f(x),g(y))$ is a solution to the both equations. 
These examples are provided by Militaru \cite[Examples 2.4]{Mi98} and in particular they lie in the class of the well-known Lyubashenko solutions \cite{Dr92}.

Finding set-theoretical solutions of the Yang-Baxter equation is equivalent to determining \emph{set-theoretical solutions of the braid equation}, i.e., maps $r:S\times S\to S\times S$ such that the relation
\begin{align*}
r_{23}\, r_{12}\, r_{23}\, =\, r_{12}\, r_{23}\, r_{13}
\end{align*}
holds. 
In particular, a map $s$ is a set-theoretical solution of the quantum Yang-Baxter equation if and only if $r=\tau s$ is a set-theoretical solution of the braid equation.
This map $r$ is usually written as $r(x,y)=(\lambda_x(y), \rho_y(x))$ with $\lambda_x, \rho_y$ maps from $S$ into itself. Since the late 1990’s a large number of works related to this equation has been produced, including the seminal papers of Gateva-Ivanova and Van den Bergh \cite{GaVa98}, Etingof, Schedler, and Soloviev \cite{ESS99}, and Lu, Yan, and Zhu \cite{LuYZ00}. In particular, the class of \emph{involutive non-degenerate solutions}, i.e., $r^2=\id$ and the maps $\lambda_x, \rho_x$ are bijective for every $x \in S$, has been the most studied. Some algebraic structures related to the braid equation have been introduced and investigated over the years.  
Just to name a few, we remind those introduced by Rump, that are cycle sets \cite{Ru05}
and a generalization of radical rings, the braces \cite{Ru07}.
The former are in bijective correspondence with involutive non-degenerate solutions and they are also studied in recent works such as \cite{CaCaMiPi19x, CaCaPi18, CaCaPi19, LeVe17}. 
The latter are a useful tool to obtain solutions and they are widely investigated in \cite{CeGaSm18, CeJeOk14, Ga18, Sm18}.\\
Recently, generalizations of braces afford to determine bijective solutions that are not necessarily involutive, for instance through skew braces \cite{GuVe17, SmVe18}, or not bijective solutions that are left non-degenerate, such as by means of semi-braces, see \cite{CaCoSt17} and \cite{JeAr19}. In this way, Guarnieri and Vendramin \cite{GuVe17} determine solutions that arise from finite skew brace such that $r^n = \id$. Under mild assumptions, Catino, Colazzo, and Stefanelli \cite{CaCoSt19} prove that the solution $r$ associated to a semi-brace satisfy $r^{n} = r$ with $n$ natural number closely linked with the semi-brace structure. \\
During the last years, Lebed \cite{Le17} and Matsumoto and Shimizu \cite{MaSh18} deal with the particular class of \emph{idempotent solutions}, i.e., $r^2=r$, which are useful tools in the study of some algebraic structures, such as factorizable monoids. 
Furthermore, Cvetko-Vah and Verwimp \cite{Charl19} provide solutions $r$ by means of skew lattices that are \emph{cubic solutions}, i.e., $r^3 = r$. Moreover, Jespers and Van Antwerpen obtain examples of such solutions through specific semi-braces, see \cite[Theorem 5.1]{JeAr19}. In addition, examples of solutions $r$ with the property $r^4=\id$ are given by Yang in \cite[p. 511]{Ya16}.

The aim of this work is to determine new solutions of the Yang-Baxter equation by means of solutions to the pentagon equation. Specifically, we show how to involve two specific solutions of the pentagon equation defined on semigroups $S$ and $T$ in order to construct new solutions of the braid equation on the cartesian product $S\times T$.\\ In addition to this, among the solutions of the pentagon equation we characterize the special class of those that are also solutions of the quantum Yang-Baxter equation. 
Our characterization leads to consider specific semigroups belonging to the variety of semigroups $\mathcal{S}=[xyz=xwyz]$ that one can deepen in \cite{Mo03}. 
The solutions defined on such semigroups are different from those known until now. Namely, such solutions $r$ satisfy the property $r^5=r^3$ and the powers of $r$ are still solutions. 

In view of all this, we also introduce a construction of solutions of the pentagon equation on semigroups. Specifically, if $s$ and $t$ are two given solutions on semigroups $S$ and $T$, respectively, we provide suitable conditions to obtain a new solution $s\bowtie t$ on the cartesian product $S\times T$, named the \emph{matched solution} of $s$ and $t$.
In particular, these solutions are defined on the \emph{matched semigroup} $S\bowtie T$, that includes the classical Zappa product in \cite{Ku83}.

Finally, the last section is devoted to some remarks and questions about solutions of the Yang-Baxter equation that arise throughout in this work.

\section{New solutions of the pentagon equation} \label{Sezione1}

This section is devoted to introducing a new construction of solutions of the pentagon equation defined on a particular semigroup that has the cartesian product of two semigroups as underlying set. In this way, examples of such solutions can be obtained starting from classical Zappa products of semigroups.
\medskip

At first, we remind some available results in literature of the pentagon equation and some examples and constructions of solutions of this equation.\\
Hereinafter, we call any set-theoretical solution of the pentagon equation simply a \emph{PE solution}. According to the notation introduced in \cite[Proposition 8]{CaMaMi19x}, given a set $S$ and a map $s$ from $S \times S$ into itself, we write
\begin{center}
	$s(x,y)=(x y, \,\theta_x(y))$,
\end{center}
where $\theta_x$ is a map from $S$ into itself, for every $x \in S$. 
Then, $s$ is a PE solution on $S$ if and only if the following conditions hold
	\begin{align}
	(x y) z&=x  (y  z) \label{assoc}\\
	\theta_x(y)  \theta_{x  y}(z)&=\theta_x(y z)\label{omotheta}\\
	\theta_{\theta_x(y)}\theta_{x  y}&=\theta_y\label{idemtheta}
	\end{align}
	for all $x,y,z \in S$. Moreover, the map $s$ is a PE solution if and only if $t:=\tau s\tau$ satisfies
\begin{equation*}
t_{12}\, t_{13}\, t_{23}=t_{23}\, t_{12}.
\end{equation*}
According to \cite[Definition 1]{CaMaMi19x}, such a map $t$ is called a \emph{set-theoretical solution of the reversed pentagon equation}, or briefly a \emph{reversed solution}.
 Thus, every reversed solution will be written as $t(x,y)=(\theta_y(x), yx)$.

\medskip

The following are easy examples of PE solutions that we will often use in this work.
\begin{exs}\label{exsKey}\hspace{1mm}
	\begin{enumerate}
		\item Let $S$ be a semigroup and $\gamma$ an idempotent endomorphism of $S$. Then, the map $s:S \times S \rightarrow S \times S$ given by 
		$$ s(x,y)=\left(x y,\gamma\left(y\right)\right)$$ 
		is a PE solution on $S$. In particular, if $e$ is an idempotent element of $S$, then $s(x,y)=(xy,e)$ is a solution on $S$.
		\item[2.] If $S$ is a set and $f,g$ are idempotent maps from $S$ into itself such that $fg=gf$, then the map $s:S \times S \rightarrow S \times S$ given by 
		$$s(x,y)=\left(f\left(x\right),\, g\left(y\right)\right)$$
		is both a PE solution and a reversed solution on $S$ that we call \emph{Militaru solution} \cite{Mi98}.
	\end{enumerate}
\end{exs}

\medskip

In literature, there are few systematic constructions of PE solutions. In \cite[Proposition 1]{KaSe98} Kashaev and Sergeev provide a construction of PE solutions on a closed under multiplication subset $S$ of a group $G$. Recently, a complete description of PE solutions $s$ on groups $G$ of the form $s(x,y)=(xy, \theta_x(y))$ is presented in \cite[Theorem 15]{CaMaMi19x}. 

The property \eqref{assoc} suggests to look for constructions of PE solutions starting from fixed semigroups and then to find maps $\theta_x$ satisfying properties \eqref{omotheta} and \eqref{idemtheta}.
\medskip

In conformity with this idea, we show a new method to obtain PE solutions on inflations of semigroups.
In detail, if $X$ is a set, $T$ a semigroup, and $\varphi: X \rightarrow T$ a map, consider $S=X \;\dot{\cup} \;T$ and $\bar{\varphi}:S \to T$ the extension map of $\varphi$ such that $\bar{\varphi}_{|_T}=\id_T$.
Thus, the set $S$ endowed with the operation given by $$a b:=\bar{\varphi}(a)\bar{\varphi}(b)$$
is a semigroup that is called the \emph{inflation of $T$ via $\varphi$}. Note that this definition is equivalent to that provided by Clifford and Preston in \cite[p. 98]{ClPr61}.\\
Given an inflation $S$ of a semigroup $T$ via a map $\varphi$, if $s(u,v)=(uv, \theta_u(v))$ is a PE solution on $T$, then the map $\bar{s}: S \times S \to S \times S$ defined by 
	\begin{equation*}
	\bar{s}\left(a,b\right)=\left(ab,  \theta_{\bar{\varphi}\left(a\right)}\bar{\varphi}\left(b\right)\right)
	\end{equation*}
	is a PE solution on the inflation $S$ of $T$. We name such a solution $\bar{s}$ an \emph{inflation of $s$ via $\varphi$}. 
\medskip

Now, we provide a more elaborate construction of PE solutions on the cartesian product of two semigroups. At first we need to fix some notations and to give preparatory definitions. For the ease of the reader, given two semigroups $S$ and $T$, we use the letters $a,b,c$ for $S$ and $u,v,w$ for $T$. Moreover, we denote any PE solution $s$ on $S$ and any PE solution $t$ on $T$ by $s(a,b)=(ab, \theta_a(b))$ and $t(u,v)=(uv,\theta_u(v))$, respectively.
\begin{defin}
	Let $S$ and $T$ be semigroups, $\alpha:T\to S^{S}$ and $\beta:S\to T^{T}$ maps, set $\alpha_{u}:= \alpha\left(u\right)$, for every $u\in T$, and $\beta_{a}:= \beta\left(a\right)$, for every $a\in S$. If $\alpha$ and $\beta$ satisfy the following conditions
	\begin{align}
	\alpha_{u}\left(a\alpha_{v}\left(b\right)\right)
	&= \alpha_{u}\left(a\right)\alpha_{\beta_{a}\left(u\right)v}\left(b\right)\label{S1} \tag{S1} \\
	\beta_{a}\left(\beta_{b}\left(u\right)v\right)
	&= \beta_{b\alpha_{v}\left(a\right)}\left(u\right)\beta_{a}\left(v\right) \tag{S2} \label{S2},
	\end{align}
	for all $a,b\in S$ and $u,v\in T$, then we call $\left(S,T,\alpha,\beta\right)$ a \emph{matched quadruple of semigroups}.
\end{defin}

It is a routine computation to verify that $S \times T$ is a semigroup with respect to the operation defined by
	\begin{equation*}
	\left(a,u\right)\left(b,v\right)=\left(a\alphaa{u}{\left(b\right)}, \betaa{b}{\left(u\right)}v\right)
	\end{equation*}
	if and only if $\left(S,T,\alpha,\beta\right)$ is a matched quadruple of semigroups. 
	We call such a semigroup the \emph{matched product of $S$ and $T$} and we denote it by $S \bowtie T$.\\
	As a class of examples of matched product of semigroups one can easily find the classical \emph{Zappa product} \cite{Ku83}. Specifically, instead of \eqref{S1} and \eqref{S2}, in this case we require the following conditions
				\begin{align}
				&\alpha_{u}\left(ab\right)
	= \alpha_{u}\left(a\right)\alpha_{\beta_{a}\left(u\right)}\left(b\right) 
	&\alpha_{uv}=\alpha_u\alpha_v \tag{S1'}\\
	&\betaa{a}{\left( uv \right)}=\betaa{\alphaa{v}{\left(a\right)}}{\left( u \right)}\betaa{a}{\left( v\right)}  &\beta_{ab}=\beta_b\beta_a \tag{S2'}
				\end{align}
hold, for all $a,b \in S$ and $u,v \in T$. 

\medskip

\begin{defin}\label{mq-sol}
	Let $\left(S,T,\alpha,\beta\right)$ be a matched quadruple of semigroups, $s$ and $t$ PE solutions on $S$ and $T$, respectively. If $s,t,\alpha$, and $\beta$ satisfy the following conditions 
	\begin{align}
	&\thetaa{a}{\alphaa{u}{}} = \thetaa{\alphaa{v}{\left(a\right)}}{\alphaa{\betaa{a}{\left(v\right) u}}{}} \tag{M1} \label{P1}\\ &\thetaa{a\alphaa{u}{\left(b\right)}}{}=\alphaa{\thetaa{\betaa{b}{\left(u\right)}}{\left(v\right)}}{\thetaa{a\alphaa{u}{\left(b\right)}}{}} \tag{M2} \label{P2}\\
	&\betaa{\thetaa{a\alphaa{u}{\left(b\right)}}{\alphaa{\betaa{b}{\left(u\right)v}}{\left(c\right)}}}{\thetaa{\betaa{b}{\left(u\right)}}{\left(v\right)}}=\thetaa{\betaa{b\alphaa{v}{\left(c\right)}}{\left(u\right)}}{\betaa{c}{\left(v\right)}} \tag{M3} \label{P3},
	\end{align}
	for all $a,b,c\in S$ and $u,v\in T$, then we call $\left(s,t,\alpha,\beta\right)$ a \emph{matched quadruple}.
\end{defin}

Now, we show how to obtain new PE solutions by means of a matched quadruple. For simplicity, we denote $S \times T \times S \times T$ by $(S\times T)^2$.

\begin{theor}\label{TeoMatch}
	Let $S$, $T$ be semigroups and $\left(s,t,\alpha,\beta\right)$ a matched quadruple. Then, the map $s \bowtie t:(S\times T)^2\to (S\times T)^2$ defined by
	\begin{align*}
	s \bowtie t\left(a,u\, ; \, b,v\right)
	= \left(a\alphaa{u}{\left(b\right)}, \betaa{b}{\left(u\right)}v\, ;\,  \thetaa{a}{\alphaa{u}{\left(b\right)}}, \thetaa{\betaa{b}{\left(u\right)}}{\left(v\right)}\right), 
	\end{align*}
	for all	$\left(a,u\right),\left(b,v\right)\in S\times T$, is a PE solution on the matched product $S \bowtie T$. We call $s\bowtie t$ the \emph{matched product of $s$ and $t$}.
	\begin{proof}
		First of all, the condition \eqref{assoc} is straightforward since $S \bowtie T$ is a semigroup.\\
	Now, if $\left(a,u\right), \left(b,v\right), \left(c,w\right) \in S \times T$, we have
		\begin{align*}
		&\thetaa{\left(a, u\right)}{\left(b, v \right)}\thetaa{\left(a,u \right)\left(b,v \right)}{\left(c,w \right)}\\
		&=\thetaa{\left(a, u\right)}{\left(b, v \right)}\thetaa{\left(a\alphaa{u}{\left(b\right)}, \betaa{b}{\left(u\right)}v\right)}{\left(c,w \right)}\\
		&=\left(\thetaa{a}{\alphaa{u}{\left(b\right)}}, \thetaa{\betaa{b}{\left(u\right)}}{v}\right)\left(\thetaa{a\alphaa{u}{\left(b\right)}}{\alphaa{\betaa{b}{\left(u\right)}v}{\left(c\right)}}, \thetaa{\betaa{c}{\left(\betaa{b}{\left(u\right)}v\right)}}{\left( w\right)}\right)\\
		&=\left(\thetaa{a}{\alphaa{u}{\left(b\right)}} \alphaa{\thetaa{\betaa{b}{\left(u \right)}}{\left(v\right)}}{\thetaa{a\alphaa{u}{\left( b \right)}}{\alphaa{\betaa{b}{\left(u \right)}v}{\left(c \right)}}}, \betaa{\thetaa{a\alphaa{u}{\left(b \right)}}{\alphaa{\betaa{b}{\left( u \right)}v}{\left( c \right)}}}{\thetaa{\betaa{b}{\left( u \right)}}{\left(v \right)}} \thetaa{\betaa{c}{\left(\betaa{b}{\left(u\right)}v\right)}}{\left( w\right)}\right)
		\end{align*}
		and
		\begin{align*}
		&\thetaa{\left(a, u \right)}{\left(\left( b,v\right)\left(c,w\right)\right)}\\
		&=\thetaa{\left(a, u \right)} {\left(b\alphaa{v}{\left(c\right)}, \betaa{c}{\left(v\right)}w\right)}\\
		&=\left( \thetaa{a}{\alphaa{u}{\left( b\alphaa{v}{\left( c\right)}\right)}} ,  \thetaa{\betaa{b \alphaa{v}{\left(c\right)}}{\left(u\right)}}{\left(\betaa{c}{\left(v\right)}w\right)} \right).
		\end{align*}
		In addition, 
		\begin{align*}
		&\thetaa{a}{\alphaa{u}{\left( b\alphaa{v}{\left( c\right)}\right)}}\\
		&=\thetaa{a}{\left(\alphaa{u}{\left(b\right)}\alphaa{\betaa{b}{\left(u\right)}v}{\left(c \right)}\right)}&\mbox{by \eqref{S1}}\\
		&=\thetaa{a}{\alphaa{u}{\left(b\right)}}\thetaa{a\alphaa{u}{\left(b\right)}}{\left(\alphaa{\betaa{b}{\left(u\right)}v}{\left(c \right)}\right)} & \\
		&=\thetaa{a}{\alphaa{u}{\left(b\right)}} \alphaa{\thetaa{\betaa{b}{\left(u \right)}}{\left(v\right)}}{\thetaa{a\alphaa{u}{\left( b \right)}}{\alphaa{\betaa{b}{\left(u \right)}v}{\left(c \right)}}} &\mbox{by \eqref{P2}}
		\end{align*}
		and
		\begin{align*}
		&\thetaa{\betaa{b \alphaa{v}{\left(c\right)}}{\left(u\right)}}{\left(\betaa{c}{\left(v\right)}w\right)} \\
		&= \thetaa{\betaa{b \alphaa{v}{\left(c\right)}}{\left(u\right)}}{\betaa{c}{\left(v\right)} }\thetaa{\betaa{b \alphaa{v}{\left(c\right)}}{\left(u\right)}\betaa{c}{\left(v\right)}}{\left( w\right)}\\
		&=\thetaa{\betaa{b \alphaa{v}{\left(c\right)}}{\left(u\right)}}{ \betaa{c}{\left(v\right)}}\thetaa{\betaa{c}{\left(\betaa{b}{\left( u \right)}v \right)}}{\left( w\right)} &\mbox{by \eqref{S2}}\\
		&=\betaa{\thetaa{a\alphaa{u}{\left(b \right)}}{\alphaa{\betaa{b}{\left( u \right)}v}{\left( c \right)}}}{\thetaa{\betaa{b}{\left( u \right)}}{\left(v \right)}} \thetaa{\betaa{c}{\left(\betaa{b}{\left(u\right)}v\right)}}{\left( w\right)} &\mbox{by \eqref{P3}}
		\end{align*}
		hence the condition \eqref{omotheta} is satisfied. Furthermore, we have
		\begin{align*}
		&\thetaa{\thetaa{\left(a, u\right)}{\left(b, v\right)}}{\thetaa{\left(a,u\right)\left(b,v\right)}{\left(c,w\right)}}\\
		&=\thetaa{\thetaa{\left(a, u\right)}{\left(b, v\right)}}{\thetaa{\left( a\alphaa{u}{\left( b\right)},\betaa{b}{\left(u\right)}v  \right)}  {\left(c,w\right)}}\\
		&= \thetaa{\left(\thetaa{a}{\alphaa{u}{\left( b \right)}},\thetaa{\betaa{b}{\left(u\right)}}{\left(v \right)}\right)}{\left(  \thetaa{a\alphaa{u}{\left(b \right)}}{\alphaa{\betaa{b}{\left( u \right)}v}{\left(c \right)}} , \thetaa{\betaa{c}{\left(\betaa{b}{\left( u\right)}v\right)}}{\left( w\right)}  \right)}\\
		&=\left(\thetaa{\thetaa{a}{\alphaa{u}{\left( b\right)}}}{\alphaa{\thetaa{\betaa{b}{\left( u\right)}}{\left(v \right)}}{\thetaa{a\alphaa{u}{\left( b\right)}}{\alphaa{\betaa{b}{\left(u\right)}v}{\left( c\right)}}}} ,\thetaa{\betaa{\thetaa{a\alphaa{u}{\left(b \right)}}{\alphaa{\betaa{b}{\left(u\right)}v}{\left(c \right)} }}{\thetaa{\betaa{b}{\left( u\right)}}{\left( v\right)}}}{\thetaa{\betaa{c}{\left(\betaa{b}{\left( u \right)}v \right)}}{\left( w\right)}}\right)
		\end{align*}
		and
		\begin{align*}
		\thetaa{\left(b,v\right)}{\left(c,w\right)}=\left(\thetaa{b}{\alphaa{v}{\left(c\right)}}, \thetaa{\betaa{c}{\left( v\right)}}{\left( w\right)}\right).
		\end{align*}
		Comparing the two relations above, we get
		\begin{align*}
		&\thetaa{b}{\alphaa{v}{\left(c\right)}}\\
		&=\thetaa{\alphaa{u}{\left(b\right)}}{\alphaa{\betaa{b}{\left(u\right) v}}{\left(c\right)}} &\mbox{by \eqref{P1}}\\
		&=\thetaa{\thetaa{a}{\alphaa{u}{\left( b\right)}}}{\thetaa{a\alphaa{u}{\left( b \right)}}{\alphaa{\betaa{b}{\left( u \right)}v}{\left( c \right)}}} &\mbox{by \eqref{idemtheta}}\\
		&=\thetaa{\thetaa{a}{\alphaa{u}{\left( b\right)}}}{\alphaa{\thetaa{\betaa{b}{\left( u\right)}}{\left(v \right)}}{\thetaa{a\alphaa{u}{\left( b\right)}}{\alphaa{\betaa{b}{\left(u\right)}v}{\left( c\right)}}}} &\mbox{by \eqref{P2}}
		\end{align*}
		and
		\begin{align*}
		&\thetaa{\betaa{c}{\left( v\right)}}{\left( w\right)}\\
		&=\thetaa{\thetaa{\betaa{b\alphaa{v}{\left( c\right)}}{\left(u \right)}}{\left( \betaa{c}{\left(v\right)}\right)}}{\thetaa{\betaa{b\alphaa{v}{\left( c\right)}}{\left( u \right)}\betaa{c}{\left(v\right)}}{\left(w \right)}}&\mbox{by \eqref{idemtheta}}\\
		&=\thetaa{\thetaa{\betaa{b\alphaa{v}{\left( c\right)}}{\left(u \right)}}{\left( \betaa{c}{\left(v\right)}\right)}}{\thetaa{\beta_{c}\left(\beta_{b}\left(u\right)v\right)}{\left(w \right)}} &\mbox{by \eqref{S2}}\\
		&=\thetaa{\betaa{\thetaa{a\alphaa{u}{\left(b \right)}}{\alphaa{\betaa{b}{\left(u\right)}v}{\left(c \right)} }}{\thetaa{\betaa{b}{\left( u\right)}}{\left( v\right)}}}{\thetaa{\betaa{c}{\left(\betaa{b}{\left( u \right)}v \right)}}{\left( w\right)}} &\mbox{by \eqref{P3}}
		\end{align*}
	thus the condition \eqref{idemtheta} holds. Therefore, the map $s\bowtie t$ is a PE solution on $S \bowtie T$.
	\end{proof}
\end{theor}

\vspace{3mm}

\begin{ex}\label{ex-match-quadr}
   Let $S,T$ be semigroups, $s$ a PE solution on $S$, and $t$ a PE solution on $T$ both in the class of solutions in \cref{exsKey}.1 given by $s(a,b)=(ab, \gamma(b))$ and  $t(u,v)=(uv,\delta(v))$, respectively. Assuming $\alpha_u=\gamma$ and $\beta_a=\delta$, for all $u \in T$ and $a \in S$, we obtain that $(s,t,\alpha, \beta)$ is a matched quadruple. Then, by \cref{TeoMatch} the map
\begin{align*}
	s \bowtie t\left(a,u\, ; \, b,v\right)= \left(a\gamma\left(b\right), \delta(u)v\, ; \, \gamma\left(b\right), \delta(v)\right)
	\end{align*}
is the matched product of $s$ and $t$.
\end{ex}

The following PE solution is defined on a matched product $S\bowtie T$ that do not lie in the class of Zappa product of semigroups.

\begin{ex}
Let $f$ be an idempotent map from $\mathbb{N}_{0}$ into itself, $S$ the semigroup on the positive integers $\mathbb{N}_{0}$ with multiplication defined by $ab = f(a)$, and $T$ the monoid $\left(\mathbb{N}_{0}, +\right)$. 
    If $\alpha_{u}=f$, for every $u \in T$, and $\beta_{a}:T\to T$ is the map defined by $\beta_{a}\left(u\right) = f\left(a\right) + u$, then  $(S,T,\alpha,\beta)$ is a matched quadruple of semigroups. 
    Indeed, if $a,b,u,v\in\mathbb{N}_{0}$, we have that $\alpha_{u}\left(a\alpha_{v}\left(b\right)\right)
 	= f^2(a)
	= \alpha_{u}\left(a\right)\alpha_{\beta_{a}\left(u\right)v}\left(b\right)$,
	hence \eqref{S1} holds, and 
    \begin{align*}
    \beta_{a}\left(\beta_{b}\left(u\right) + v\right) 
	    &= f\left(a\right) + f\left(b\right) + u + v
	    = f^2\left(b\right) + u + f\left(a\right) + v\\
	    &= \beta_{f(b)}\left(u\right) + \beta_{a}\left(v\right)
	    = \beta_{b\alpha_{v}\left(a\right)}\left(u\right) + \beta_{a}\left(v\right).
	\end{align*}
	Moreover, if $s$ is the PE solution on $S$ defined by $s(a,b)=(f(a),f(b))$ and  $t$ the PE solution on $T$ given by $t(u,v)=(uv,v)$ it follows that $(s,t,\alpha,\beta)$ is a matched quadruple. Indeed, \eqref{P1} and \eqref{P2} are trivially satisfied. Furthermore, it holds
    \begin{align*}
        \betaa{\thetaa{a\alphaa{u}{\left(b\right)}}{\alphaa{\betaa{b}{\left(u\right) + v}}{\left(c\right)}}}{\thetaa{\betaa{b}{\left(u\right)}}{\left(v\right)}}= \betaa{f^2\left(c\right)}{\left(v\right)}
        = \betaa{f\left(c\right)}{\left(v\right)}
        = f^{2}\left(c\right) + v
        = \betaa{c}{\left(v\right)}
        =\thetaa{\betaa{b\alphaa{v}{\left(c\right)}}{\left(u\right)}}{\betaa{c}{\left(v\right)}},
    \end{align*}
    i.e., the condition $\eqref{P3}$ is satisfied. Therefore, by \cref{TeoMatch} the map defined by
    \begin{align*}
        s\bowtie t\left(a,u; b,v\right)
        = \left(f(a), f\left(b\right) + u + v \, ; \, f\left(b\right), v\right)
    \end{align*}
    is the matched product of $s$ and $t$.
\end{ex}

\vspace{3mm}

Given a matched product of solutions $s \bowtie t$, in order to extrapolate the solutions $s$ and $t$, we need to find isomorphic copies of $S$ and $T$ inside $S\bowtie T$. For this purpose, we have to require additional properties on semigroups $S$ and $T$. 

\begin{prop}\cite[cf. Theorem 5]{Ca87}\label{th:Cat87}
	Let $S \bowtie T$ be a matched product of semigroups via $\alpha$ and $\beta$, $e_S$ a right identity for $S$, and $e_T$ a left identity for $T$. Then, if the following conditions
	\begin{align}
	&\alpha_{e_T}=\id_S \label{S3} \tag{S3}\\
	&\forall \ a \in S \quad \betaa{a}{\left(e_T\right)}=e_T \label{S4} \tag{S4}
	\end{align}
	hold, then $S^*=\left\lbrace \left(a,e_T\right) \; \vert \; a \in S\right\rbrace$ is a subsemigroup of $S \bowtie T$ isomorphic to $S$. \\
	Similarly, if
	\begin{align}
	&\beta_{e_S}=\id_T \label{S5} \tag{S5}\\
	&\forall \ u \in T \quad \alphaa{u}{\left( e_S \right)}=e_S \label{S6} \tag{S6}
	\end{align}
	hold, then $T^*=\left\lbrace \left(e_S,u\right) \; \vert \; u \in T\right\rbrace$ is a semigroup of $S \bowtie T$ isomorphic to $T$.
\end{prop}

\begin{ex} Let $X=\lbrace 0,1,2 \rbrace$, $T$ the right zero semigroup on $X$, and $S$ the semigroup on $X$ with the multiplication given by
\begin{align*}
  &\forall\, a \in X \quad  a \cdot 0=0 \cdot a=0,\\
  &\forall\, a,b \in X \setminus \lbrace 0 \rbrace \quad a\cdot  b=a.
\end{align*}
Then, $e_S=1$ is a right identity for $S$ and $e_T=0$ is a left identity for $T$. It is a routine computation to check that the map $\gamma: X \to X$ defined by $\gamma(x)=1\cdot  x$ is an idempotent endomorphism of $S$. Moreover, let $\alpha: T \to S^S$ and $\beta: S \to T^T$ be the maps such that
\begin{center}
		$\alpha_u=\begin{cases}
		\id_S \qquad &\text{if} \;u=0 \\
		\gamma \quad&\text{if}\; u \neq 0
		\end{cases}$ \qquad $\beta_a=\begin{cases}
		\id_T \qquad &\text{if} \;a=1 \\
		\gamma \quad&\text{if}\; a \neq 1
		\end{cases}$
	\end{center}
	Then, $S \bowtie T$ is a matched product of semigroups via the maps $\alpha$ and $\beta$ (that do not lie in the class of Zappa product). Furthermore, conditions \eqref{S3}, \eqref{S4}, \eqref{S5}, and \eqref{S6} are trivially satisfied and thus, by \cref{th:Cat87}, one can find isomorphic copies of $S$ and $T$ inside $S \bowtie T$. \\
	Let $s$ be the PE solution on $S$ given by $s(a,b)=(a\cdot b, \gamma(b))$ and $t$ the PE solution on $T$ given by $t(u,v)=(v,v)$. Then, it is a routine computation to check that \eqref{P1}, \eqref{P2}, and \eqref{P3} and so by \cref{TeoMatch} the map
	\begin{align*}
	s \bowtie t\left(a,u\, ; \, b,v\right)= \left(a\cdot\alpha_u\left(b\right), v\, ; \, \gamma\left(b\right), v\right)
	\end{align*}
	is the matched product of $s$ and $t$. 
\end{ex}

\begin{rem} \label{rem-semi} Let $(S,T,\alpha,\beta)$ be a matched quadruple of semigroups and $e_S$ a right identity for $S$. If conditions \eqref{S5} and \eqref{S6} hold, then the condition \eqref{P2} becomes easier and it is equivalent to
	\begin{align}
	\theta_a=\alphaa{\theta_u(v)}{\theta_a},\tag{M2$'$} \label{P2'}
	\end{align}
	for all $a \in S$ and $u,v \in T$. Indeed, if \eqref{P2} holds, with $b=e_S$, we get
	\begin{align*}
	\theta_a=\theta_{a e_S}\underset{\eqref{S6}}{=}\theta_{a\alphaa{u}{(e_S)}}=\alphaa{\thetaa{\betaa{e_S}{(u)}}{(v)}}{\thetaa{a\alphaa{u}{(e_S)}}{}}\underset{\eqref{S5},\eqref{S6}}{=}\alphaa{\thetaa{u}{(v)}}{\thetaa{ae_S}{}}=\alphaa{\thetaa{u}{(v)}}{\thetaa{a}{}},
	\end{align*}
	i.e., \eqref{P2'} holds. Conversely, if \eqref{P2'} holds, with $a=a\alphaa{u}{(b)}$ and $u=\betaa{b}{(u)}$, we clearly obtain that condition \eqref{P2} is satisfied.
\end{rem}

If $S$ and $T$ are monoids and the identities $1_S$ and $1_T$ satisfy the conditions \eqref{S3}, \eqref{S4}, \eqref{S5}, \eqref{S6}, we call any matched quadruple $\left(S,T,\alpha,\beta\right)$ a \emph{matched quadruple of monoids}. Moreover, $S \bowtie T$ is a monoid with identity $\left(1_S, 1_T\right)$ if and only if $(S,T,\alpha,\beta)$ is a matched quadruple of monoids. Note that in such a case this construction is equivalent to the classical Zappa product of two monoids.

In the following we show that under the assumption of $(S,T,\alpha, \beta)$ is a matched quadruple of monoids, then the conditions in \cref{mq-sol} become easier. 
\begin{prop}\label{prop-match-mon}
	Let $(S,T,\alpha,\beta)$ be a matched quadruple of monoids. Then, $(s,t,\alpha,\beta)$ is a matched quadruple if and only if the following conditions
	\begin{enumerate}
		\item $\theta_a=\alphaa{\theta_u(v)}{\theta_a}=\thetaa{\alphaa{v}{(a)}}{\alphaa{\betaa{a}{(v)}}{}}$,
		\item $\betaa{\thetaa{a}{\alphaa{uv}{(b)}}}{\thetaa{u}{(v)}}=\thetaa{\betaa{\alphaa{v}{(b)}}{(u)}}{\betaa{b}{(v)}}$,
	\end{enumerate}
	hold, for all $a, b \in S$ and $u, v\in T$.
	\begin{proof}
		Initially, note that by \cref{rem-semi} the first equality of condition $1.$ is equivalent to the condition \eqref{P2}.\\
		Now, suppose that $(s,t,\alpha,\beta)$ is a matched quadruple. If $a \in S$ and $v \in T$, then
		\begin{align*}
		\thetaa{a}{}
		    &= \theta_a\alpha_{1_T}
		    &\mbox{by \eqref{S3}}\\
		    &= \thetaa{\alphaa{v}{(a)}}{\alphaa{\betaa{a}{(v)}1_T}{}}
		    &\mbox{by \eqref{P1}}\\
		    &=\thetaa{\alphaa{v}{(a)}}{\alphaa{\betaa{a}{(v)}}{}}.
		\end{align*}
		Moreover, we have
		\begin{align*}
		&\betaa{\thetaa{a}{\alphaa{uv}{(b)}}}{\thetaa{u}{(v)}}\\
		&= \betaa{\thetaa{a\alphaa{u}{\left(1_S\right)}}{\alphaa{\betaa{1_S}{\left(u\right)v}}{\left(b\right)}}}{\thetaa{\betaa{1_S}{\left(u\right)}}{\left(v\right)}} &\mbox{by \eqref{S5}-\eqref{S6}}\\
		&= \thetaa{\betaa{1_S\alphaa{v}{\left(b\right)}}{\left(u\right)}}{\betaa{b}{\left(v\right)}}&\mbox{by \eqref{P3}}\\
		&=\thetaa{\betaa{\alphaa{v}{(b)}}{(u)}}{\betaa{b}{(v)}}.
		\end{align*}
		Conversely, suppose that conditions $1.$ and $2.$ hold. Then, if $a, b \in S$ and $u,v \in T$, since $\alpha$ is a homomorphism, we get
		\begin{align*}
		\thetaa{a}{\alphaa{u}{}} = \thetaa{\alphaa{v}{(a)}}{\alphaa{\betaa{a}{(v)}}{}}\alpha_u= \thetaa{\alphaa{v}{\left(a\right)}}{\alphaa{\betaa{a}{\left(v\right) u}}{}},
		\end{align*}
		i.e., the condition \eqref{P1} holds. Furthermore, set $a=a\alphaa{u}{(b)}$, $u=\betaa{b}{(u)}$, $c=b$, we trivially obtain the condition \eqref{P3}.
	\end{proof}
\end{prop}

\begin{ex}
Let $S=\lbrace 1_S,x,y\rbrace$ be the monoid such that $x^2=x, \, y^2=y, \,xy=y=yx$ and $T=\lbrace 1_T, z\rbrace$ the monoid such that $z^2=z$. Consider the maps $\alpha: T \to S^S$ such that $\alpha_{1_T}=\id_S$ and $\alpha_z=\gamma$, where $\gamma$ is an idempotent endomorphism of $S$, and $\beta:S\to T^T$ such that $\beta_{1_S}=\id_T$ and $\beta_x(u)=\beta_y(u)=1_T$, for every $u \in T$. Then, it is easy to check that $(S,T,\alpha,\beta)$ is a matched quadruple of monoids.\\
Moreover, let $s$ be the PE solution on $S$ given by $s(a,b)=(ab, \gamma(b))$ and $t$ the PE solution on $T$ given by $t(u,v)=(uv,v)$. Then, by \cref{prop-match-mon} $(s,t,\alpha, \beta)$ is a matched quadruple and so the map
    \begin{center}
		$s \bowtie t(a,u;\,b,v)=\begin{cases}
		(a,v; \, 1_S,v) \qquad &\text{if} \;b=1_S, \,u=1_T \\
		(a,uv; \, 1_S,v) \qquad &\text{if} \;b= 1_S, \,u\neq 1_T\\
		(ab,v; \, \gamma(b),v) \qquad &\text{if} \;b\neq 1_S, \,u=1_T\\
		(a\gamma(b),v; \, \gamma(b),v) \qquad &\text{if} \;b\neq 1_S, \,u\neq 1_T
		\end{cases}$ 
	\end{center}
is the matched product of $s$ and $t$.
\end{ex}

\section{QYBE solutions of pentagonal type}

This section is devoted to the special class of PE solutions on semigroups of the form $s(a,b)=(ab, \theta_a(b))$ that are also solutions of the quantum Yang-Baxter equation.  In particular, we focus on semigroups belonging to the variety $\mathcal{S}=[abc=adbc]$ (cf. \cite{Mo03}) on which we are able to describe all such solutions.

\medskip

Hereinafter, we call any set-theoretical solution of the quantum Yang-Baxter equation briefly a \emph{QYBE solution}. Initially, we provide a characterization of PE solutions of the form $s(a,b)=(ab, \theta_a(b))$ that are also QYBE solutions.

\begin{prop}\label{prop-PE-YBE}
	Let $S$ be a semigroup and $s$ a PE solution on $S$ defined by $s(a,b)=(ab,\theta_a(b))$. Then, the map $s$ is a QYBE solution if and only if the following conditions	
	\begin{align}
	&\label{YB1}\tag{Y1}abc=a\theta_b(c)bc\\
	&\label{YB2}\tag{Y2}\theta_a\theta_b = \theta_b\\	
	&\label{YB3}\tag{Y3}\thetaa{a}{(bc)}=\thetaa{\theta_b(c)}{(bc)}
	\end{align}
	are satisfied, for all $a,b,c \in S$.
	\begin{proof}
 If $a,b,c \in S$, we have
 \begin{align*}
     s_{12}s_{13}s_{23}\left(a,b,c\right)=\left(a\theta_b(c)bc, \thetaa{a\theta_b(c)}{(bc)},\theta_a\theta_b\right)
 \end{align*}
 and by \eqref{omotheta} and \eqref{idemtheta} 
 \begin{align*}
     s_{23}s_{13}s_{12}\left(a,b,c\right)=\left( abc, \theta_a\left(b\right)\theta_{ab}\left(c\right), \theta_{\theta_a\left(b\right)}\theta_{ab}\left(c\right)\right)=\left(abc, \theta_a\left(bc\right), \theta_b\left(c\right)\right).
 \end{align*}
 Assuming $s$ is a QYBE solution, we note that comparing the first and the third components, the conditions \eqref{YB1} and \eqref{YB2} hold. Moreover, by \eqref{YB2} and \eqref{idemtheta}, we obtain $\theta_{ab}=\theta_{\theta_a\left(b\right)}\theta_{ab}=\theta_b$ and so
		\begin{align}\label{theta_ab}
		\theta_{ab}=\theta_b.
		\end{align}	
Hence, comparing the second components we have that $\thetaa{\theta_b(c)}{(bc)}=\thetaa{a\theta_b(c)}{(bc)}=\theta_a\left(bc\right)$ and then \eqref{YB3} is satisfied.\\
Conversely, if \eqref{YB1}, \eqref{YB2}, and \eqref{YB3} hold, then by using \eqref{theta_ab} we get the claim.
	\end{proof}
\end{prop}

\noindent Recalling that a map $s$ is a PE solution if and only if the map $t=\tau s \tau$ is a reversed solution, we have the following result.

\begin{cor}
Let $S$ be a semigroup and $t$ a reversed solution on $S$ defined by $t(a,b)=(\theta_b(a),ba)$. Then, the map $t$ is a QYBE solution on $S$ if and only if the conditions \eqref{YB1}, \eqref{YB2}, and \eqref{YB3}
	are satisfied. 
\end{cor}

\begin{defin}
    A PE solution $s$ satisfying \eqref{YB1}, \eqref{YB2}, and \eqref{YB3} is said to be a \emph{QYBE solution of pentagonal type}, or briefly \emph{a P-QYBE solution}. Similarly, a reversed solution satisfying \eqref{YB1}, \eqref{YB2}, and \eqref{YB3} is called a \emph{QYBE  solution of reversed pentagonal type}, or briefly \emph{a R-QYBE solution}.
\end{defin}

\noindent From now on, we will show some results for P-QYBE solutions that can be equivalently obtained for R-QYBE solutions.

The following proposition shows in which way one can find the maps $\theta_a$ in order to construct a P-QYBE solution.

\begin{prop}\label{rem-prod-direct} Let $s(a,b)=(ab, \theta_a(b))$ be a P-QYBE solution on a semigroup $S$. Then, the following hold:
\begin{enumerate}
    \item the map $\theta_a$ is idempotent, for every $a \in S$;
    \item  ${\theta_a}_{|_{S^2}}={\theta_b}_{|_{S^2}}$, for all $a,b \in S$;
    \vspace{1mm}
    \item if $S^2=S$, then $s(a,b)=(ab, \bar{\theta}(b))$, where $\bar{\theta}$ is an idempotent endomorphism of $S$.
\end{enumerate}
\begin{proof}\hspace{1mm}\\
        1. The claim follows by \eqref{YB2}.\\
    	2. The claim follows by \eqref{YB3}.\\
		3. By 2., the maps $\theta_x$ are all equal. Let $x \in S$, set $\theta_x:=\bar{\theta}$, we have that $\bar{\theta}$ is an idempotent endomorphisms of $S$. Indeed, by \eqref{idemtheta} it holds $\bar{\theta}^{2}=\bar{\theta}$ and we obtain
		\begin{align*}
		\bar{\theta}(ab)=\theta_x(ab)\underset{\eqref{omotheta}}{=}\theta_x(a)\theta_{xa}(b)=\bar{\theta}(a)\bar{\theta}(b),
		\end{align*}
	for all $a,b \in S$.

\end{proof}

\end{prop}

The following are simple examples of PE solutions of the form $s(a,b)=(ab, \gamma(b))$ in \cref{exsKey}.1 that are of P-QYBE type.

\begin{exs}\label{es1: PE-YBE-1} \hspace{1mm}
    \begin{enumerate}
		\item The solution $s\left(a,b\right)=\left(ab, e_S\right)$, with $e_S$ a left identity (or a right identity) for $S$. Note that in the particular case of $S$ a group, by \eqref{YB1} the unique P-QYBE solution $s$ on $S$ is given by $s(a,b)=(ab,1)$.
		\vspace{1mm}
		\item The map $s\left(a,b\right)=\left(ab, b\right)$, with $S$ is a left quasi-normal semigroup, i.e., $abc = acbc$, for all $a,b,c, \in S$, (for more details see \cite{Go83}).
	\end{enumerate}	
\end{exs}

In order to find more examples, we note that the condition \eqref{YB1} leads to consider special classes of semigroups. Thereby, we focus on semigroups $S$ belonging to the variety 
\begin{equation}\label{GRB}
\mathcal{S}:=[abc=adbc]
\end{equation}
 in \cite{Mo03} which immediately ensures \eqref{YB1}. In this way, one has to find maps $\theta_a$ from $S$ into itself satisfying just \eqref{YB2} and \eqref{YB3}. 
 
 \begin{ex}\label{ex7}
 Let $S \in \mathcal{S}$ and $s(a,b)=(ab, \gamma(b))$ the PE solution on $S$ in \cref{exsKey}.1. Then, $s$ is a P-QYBE solution.
 \end{ex}
 As a direct consequence of \cref{rem-prod-direct} and \cref{ex7}, we provide a complete description of P-QYBE solutions in the particular case of semigroups $S \in \mathcal{S}$ such that $S^2=S$.  

\begin{prop}\label{rect-band} Let $S\in \mathcal{S}$ 
such that $S^2=S$. Then, the unique P-QYBE solutions on $S$ are of the form
\begin{align}\label{sgamma}
  s(a,b)=(ab, \bar{\theta}(b)),  
\end{align}
with $\bar{\theta}$ an idempotent endomorphism of $S$. 
\end{prop}

The following are examples of P-QYBE solutions defined on a semigroup $S\in \mathcal{S}$ for which in general $S^2\neq S$.
 
\begin{exs} \label{es-key-2} \hspace{1mm}
\begin{enumerate}
    \item The map $s\left(a,b\right)=\left(ab, bab\right)$.
    \vspace{1mm}
	\item The Militaru solution $s(a,b)=(f(a),g(b))$, where $ab=f(a)$, for all $a,b \in S$.
	\end{enumerate}
\end{exs}

\medskip

The following proposition shows some properties related to the powers of P-QYBE solutions.
\begin{prop}\label{cor-potenza}
	Let $s$ be a P-QYBE solution on $S$. Then, for every $n \in \mathbb{N}$, $n \geq 2$,
	\begin{align*}
	s^n\left(a,b\right)=\left(\ ab\,\theta_a\left(b\right)^{n-1}, \theta_a\left(b\right) \right),
	\end{align*}
	for all $a,b\in S$. In particular, if $S$ is an idempotent semigroup, it holds that $s^3=s^2$.
	\begin{proof} We prove the claim by induction on $n$. The case $n=2$ follows from the definition of $s$. Suppose that the thesis holds for $n>2$ and so by induction hypothesis, we have
		\begin{align*}
		s^{n+1}\left(a,b\right)&=s\left(ab\,\theta_a\left(b\right)^{n-1}, \theta_a\left(b\right)\right)\\
			&=\left( ab\,\theta_a\left(b\right)^{n},\theta_{ab\theta_a\left(b\right)^{n-1}}\theta_a\left(b\right)\right)&\mbox{by $\eqref{YB2}$}\\
		&=\left(ab\,\theta_a\left(b\right)^{n}, \theta_a\left(b\right)\, \right).
		\end{align*} 
		Therefore, the statement follows.	
	\end{proof}
\end{prop}

\medskip

By using \cref{cor-potenza}, we are able to provide P-QYBE solutions whose powers are still P-QYBE solutions. 

\begin{exs}\label{ex:powers}\hspace{1mm}
	\begin{enumerate}
	\item The map $s(a,b)=(ab, bab)$ that satisfies $s^3=s^2$.
	\vspace{1mm}
		
		\item  The Militaru solution $s(a,b)=(f(a),g(b))$ for which $s^2=s$ holds.
		\item The map $s(a,b)=(ab, e_S)$ in \cref{es1: PE-YBE-1}.1. In particular, it holds
		\begin{itemize}
		    \item[-] $s^2=s$, if $e_S$ is a right identity for $S$;
		    \item[-] $s^3=s^2$, if $e_S$ is a left identity for $S$.
\end{itemize}
		\item The map $s(a,b)=(ab,b)$ in \cref{es1: PE-YBE-1}.2 for which it holds $s^3=s^2$.
	\end{enumerate}
\end{exs}

\medskip

Instead, the following is an example of P-QYBE solution whose powers are not all P-QYBE solutions.
\begin{ex}
Let $S\in \mathcal{S}$ and $s$ the P-QYBE solution on $S$ given by $s(a,b)=(ab, \gamma(b))$ in \cref{ex7}, with $\gamma \neq \id_S$. Then, $s$ satisfies $s^4=s^3$ and in this case $s^3$ is a solution P-QYBE, while $s^2$ is neither a PE solution nor a QYBE solution.
\end{ex}
\bigskip

\section{Particular classes of YBE solutions}

This section is devoted to studying the powers of the braid version $r(a,b)=(\theta_a(b),ab)$ of specific P-QYBE solutions. As observed by Yang \cite[p. 16]{Ya16}, if $r$ is a YBE solution, its $n$th power $r^n$ is not necessarily a YBE solution.
In contrast to this fact, we show that the powers of solutions to the braid equation defined on semigroups $S\in \mathcal{S}$ in \eqref{GRB} are still solutions. We underline that these maps $r$ lie in the class of degenerate solutions, unless trivial cases. 

\medskip

From now on, we call any set-theoretical solution of the braid equation briefly a \emph{YBE solution} and the braid version of any P-QYBE solution simply a \emph{P-YBE solution}.\\
In the next theorem we provide sufficient conditions so that the powers of any P-YBE solution are still solutions.

\begin{theor}\label{th-power}
	Let $S \in \mathcal{S}$ and $r(a,b)=(\theta_a(b),ab)$ a P-YBE solution on $S$. Then, it holds $$r^5=r^3$$
	and the powers $r^2,r^3,r^4$ of the map $r$ are still YBE solutions. \\
	In addition, if $S$ is idempotent, it holds $r^4=r^2$.	
	\begin{proof}
		First,  denote the map ${\theta_a}_{|_{S^2}}$ as $\bar{\theta}$, for every $a \in S$. Note that, by \cref{prop-PE-YBE} conditions \eqref{YB2} and \eqref{YB3} hold. If $a,b \in S$, we have that
		\begin{align*}
		r^2\left(a,b\right)
		&=\left( \theta_{\theta_a\left(b\right)}\left(ab\right), \theta_a\left(b\right)ab \right)=\left( \bar{\theta}\left(ab\right), \theta_a\left(b\right)ab\right)
		\end{align*}
		and
		\begin{align*}
		r^3\left(a,b\right)&=\left(\bar{\theta}\left(\theta_a\left(b\right)ab\right),  \bar{\theta}\left(ab\right)\theta_a\left(b\right)ab \right)\\
		&=\left(\bar{\theta}\left(\theta_a\left(b\right)ab\right),  \bar{\theta}\left(a\right)ab\right)&\mbox{by \eqref{omotheta}-\eqref{GRB}}\\
		&=\left(\theta_a\left(b\right)\bar{\theta}\left(ab\right),  \bar{\theta}\left(a\right)ab \right)&\mbox{by \eqref{omotheta}-\eqref{YB2}}
		\end{align*}
		By proceeding in this way, one obtains
		\begin{align*}
		r^4\left(a,b\right)=\left( \bar{\theta}\left(a\right)\bar{\theta}\left(ab\right), \theta_a\left(b\right)ab\right), \qquad r^5\left(a,b\right)=\left(\theta_a\left(b\right)\bar{\theta}\left(ab\right) , \bar{\theta}\left(a\right)ab\right).
		\end{align*}
		Consequently, it holds $r^5=r^3$.\\
		In particular, if $S$ is idempotent, comparing the first components of $r^2$ and $r^4$, we get
		\begin{align*}
		\bar{\theta}\left(ab\right)&=\bar{\theta}\left(aab\right)\\
		&=\bar{\theta}(a)\theta_a(ab) &\mbox{by $\eqref{theta_ab}$} \\
		&=\bar{\theta}(a)\bar{\theta}(ab).
		\end{align*}
	 Therefore, it holds $r^4=r^2$.\\
		Now, in order to show that $r^2$ is a YBE solution, set
		\begin{align*}
		\lambda_{a}\left(b\right):=\bar{\theta}\left(ab\right)\qquad \rho_{b}\left(a\right):=\theta_a\left(b\right)ab
		\end{align*}
		and we check that
		\begin{align*}
		&\lambda_{a}\lambda_{b}\left(c\right) = \lambda_{\lambda_{a}\left(b\right)}\lambda_{\rho_{b}\left(a\right)}\left(c\right)\\
		&\rho_{c}\rho_{b}\left(a\right)
		= \rho_{\rho_{c}\left(b\right)}\rho_{\lambda_{b}\left(c\right)}\left(a\right)\\
		&\lambda_{\rho_{\lambda_{b}\left(c\right)}\left(a\right)}\rho_{c}\left(b\right)
		= \rho_{\lambda_{\rho_{b}\left(a\right)}\left(c\right)}\lambda_{a}\left(b\right),
		\end{align*}
		for all $a,b,c\in S$.
		Then, we have
		\begin{align*}
		\lambda_{a}\lambda_{b}\left(c\right)&=\lambda_{a}\bar{\theta}\left(bc\right)=\bar{\theta}\left(a\bar{\theta}\left(bc\right)\right)\\
		&=\bar{\theta}\left(a\right)\bar{\theta}\left(bc\right)&\mbox{by $\eqref{omotheta}$-$\eqref{YB2}$}              
		\end{align*}                     
		and                                          
		\begin{align*}
	\lambda_{\lambda_{a}\left(b\right)}\lambda_{\rho_{b}\left(a\right)}\left(c\right)	&=\lambda_{\bar{\theta}\left(ab\right)}\lambda_{\theta_a\left(b\right)}\left(c\right)=\lambda_{\bar{\theta}\left(ab\right)}\bar{\theta}\left(\theta_a\left(b\right)c\right)\\
		&=\bar{\theta}\left(\bar{\theta}\left(ab\right)\bar{\theta}\left(\theta_a\left(b\right)c\right)\right)\\
		&=\bar{\theta}\left(\bar{\theta}\left(a\right)\theta_a\left(b\right)\bar{\theta}\left(c\right)\right) &\mbox{by $\eqref{omotheta}$-$\eqref{GRB}$-$\eqref{YB2}$}\\
		&=\bar{\theta}\left(a\right)\bar{\theta}\left(\theta_a\left(b\right)\bar{\theta}\left(c\right)\right)&\mbox{by \eqref{omotheta}-\eqref{YB2}}\\
		&=\bar{\theta}\left(a\right)\bar{\theta}\left(bc\right)&\mbox{by $\eqref{omotheta}$-$\eqref{YB2}$}
		\end{align*}
		Moreover, we compute
		\begin{align*}
		\rho_{c}\rho_{b}\left(a\right)&=\rho_{c}\left(\theta_a\left(b\right)ab\right)=\theta_{\theta_a\left(b\right)ab}\left(c\right)\theta_a\left(b\right)abc\\
		&=\theta_b\left(c\right)bc&\mbox{by \eqref{theta_ab}-\eqref{GRB}}
		\end{align*}
		and
		\begin{align*}
		\rho_{\rho_{c}\left(b\right)}\rho_{\lambda_{b}\left(c\right)}\left(a\right)&=\rho_{\theta_b\left(c\right)bc}\rho_{\bar{\theta}\left(bc\right)}\left(a\right)=	\rho_{\theta_b\left(c\right)bc}\,\theta_a\left(\bar{\theta}\left(bc\right)a\bar{\theta}\left(bc\right)\right)\\
		&=\rho_{\theta_b\left(c\right)bc}\, \theta_a\left(	\bar{\theta}\left(b\right)\bar{\theta}\left(bc\right)\right)&\mbox{by \eqref{omotheta}-\eqref{GRB}}\\
		&=\rho_{\theta_b\left(c\right)bc}\left(	\bar{\theta}\left(b\right)\bar{\theta}\left(bc\right)\right)&\mbox{by \eqref{YB2}}\\
		&=\theta_{\bar{\theta}\left(b\right)\bar{\theta}\left(bc\right)}\left(\theta_b\left(c\right)bc\right)\bar{\theta}\left(b\right)\bar{\theta}\left(bc\right)\theta_b\left(c\right)bc\\
		&=\theta_b\left(c\right)bc&\mbox{by $\eqref{YB2}$-$\eqref{GRB}$}\\
		\end{align*}
		Finally, we have that
		\begin{align*}
		\lambda_{\rho_{\lambda_{b}\left(c\right)}\left(a\right)}\rho_{c}\left(b\right)
		&= \lambda_{\rho_{\bar{\theta}\left(bc\right)}\left(a\right)}\left(\theta_b\left(c\right)bc\right)=\lambda_{\theta_a\left(\bar{\theta}\left(bc\right)a\bar{\theta}\left(bc\right)\right)}\left(\theta_b\left(c\right)bc\right)\\
		&=\lambda_{\bar{\theta}\left(b\right)\bar{\theta}\left(bc\right)}\left(\theta_b\left(c\right)bc\right)&\mbox{by $\eqref{YB2}$-$\eqref{GRB}$} \\
		&=\bar{\theta}\left(\bar{\theta}\left(b\right)\bar{\theta}\left(bc\right)\theta_b\left(c\right)bc\right)\\
		&=\bar{\theta}\left(\bar{\theta}\left(b\right)bc\right)&\mbox{by \eqref{GRB}}\\
		&=\bar{\theta}\left(b\right)\bar{\theta}\left(bc\right)&\mbox{by \eqref{YB2}-\eqref{omotheta}}
		\end{align*}
		and
		\begin{align*}
		\rho_{\lambda_{\rho_{b}\left(a\right)}\left(c\right)}\lambda_{a}\left(b\right)
		&=
		\rho_{\lambda_{\theta_a\left(b\right)ab}\left(c\right)}\bar{\theta}\left(ab\right)=\rho_{\bar{\theta}\left(\theta_a\left(b\right)abc\right)}\bar{\theta}\left(ab\right)\\
		&=\rho_{\theta_a\left(b\right)\bar{\theta}\left(bc\right)}\bar{\theta}\left(ab\right)	\\
		&=\theta_{\bar{\theta}\left(ab\right)}\left(\theta_a\left(b\right)\bar{\theta}\left(bc\right)\right)\bar{\theta}\left(ab\right)\theta_a\left(b\right)\bar{\theta}\left(bc\right)\\
		&=\theta_{\bar{\theta}\left(ab\right)}\left(\theta_a\left(b\right)\bar{\theta}\left(bc\right)\right)\bar{\theta}\left(bc\right)&\mbox{by \eqref{GRB}}\\
		&=\theta_a\left(b\right)\bar{\theta}\left(bc\right)\bar{\theta}\left(bc\right)&\mbox{by \eqref{YB2}-\eqref{omotheta}}\\
		&=\theta_a\left(b\right)\bar{\theta}\left(bc\right)&\mbox{by \eqref{GRB}}
		\end{align*}
		By \eqref{YB3} we obtain that
		\begin{center}
			$\theta_a\left(b\right)\bar{\theta}\left(bc\right)=\theta_a\left(b\right)\theta_b\left(bc\right)=\theta_a\left(bbc\right)=\bar{\theta}\left(bbc\right)=\bar{\theta}\left(b\right)\bar{\theta}\left(bc\right)$.
		\end{center}
		Therefore, $r^2$ is a YBE solution.\\
		With similar computations one can check that $r^3,r^4$ are YBE solutions.
	\end{proof}
\end{theor}

\begin{rem}\label{rem_left_quasi_normal}
There exist P-YBE solutions $r$ for which $r^5=r^3$, but the powers of $r$ are not solutions. The P-YBE solution $r(a,b)=(b,ab)$ defined on a left quasi-normal semigroup $S$ is such an example. We highlight that this class of semigroups strictly contains $\mathcal{S}$ in \eqref{GRB}.
\end{rem}

\noindent The following example shows that the conditions in \cref{th-power} are not necessary.

\begin{ex}
	Let $S$ be a semigroup, $e_S$ a left identity for $S$, and $r$ the P-YBE solution on $S$ defined by $r\left(a,b\right)=\left(e_S, ab\right)$. Then, it holds $r^2=r$ and so clearly $r^5=r^3$.
\end{ex}

\noindent Since solutions in \cref{ex7} and \cref{es-key-2} satisfy the hypotheses of \cref{th-power}, we list the powers of the braid version of the solutions therein that are still solutions.

\begin{exs}\label{ex:powers-2}\hspace{1mm}
	\begin{enumerate}
		\item The map $r\left(a,b\right)=\left(\gamma(b), ab\right)$ satisfies $r^5=r^3$
		and the solutions $r^2,r^3,r^4$ are
		\begin{itemize}
		    \item[-] $r^2\left(a,b\right)=\left( \gamma(ab), \gamma(b)ab\right)$
		    \item[-] $r^3\left(a,b\right)=\left( \gamma(bab), \gamma(a)ab\right)$,
		\item[-] $r^4\left(a,b\right)=\left(\gamma(a^2b), \gamma(b)ab\right)$.
		\end{itemize}

		\vspace{1mm}
		
		\item The map $r(a,b)=(bab, ab)$ satisfies $r^4=r^2$ and the solutions $r^2, r^3$ are 
		\begin{itemize}
		    \item[-] $r^2\left(a,b\right)=\left( aab, bab\right)$
		    \item[-] $r^3\left(a,b\right)=\left( bab, aab\right)$.
		\end{itemize}
		\vspace{1mm}
		\item The solution
		$r\left(a,b\right)=\left(g\left(b\right), f\left(a\right)\right)$
satisfies $r^4=r^2$ and $r^2,r^3$ are respectively
		\begin{itemize}
		\item[-]$r^2\left(a,b\right)=\left( fg(a), fg(b)\right)$,
		\item[-]$r^3\left(a,b\right)=\left( fg(b), fg(a)\right)$.
		\end{itemize}
		Since $ab=f(a)$, note that in general $S$ is not idempotent and so the vice versa of the last claim of \cref{th-power} does not hold.
	\end{enumerate} 
\end{exs}

To conclude this section, we give an application of \cref{th-power} in the case of a rectangular band $S$, i.e., an idempotent semigroup such that $abc=ac$, for all $a,b,c \in S$ (cf. \cite{Ho95}).

\begin{ex} Let $S$ be a rectangular band. Then, by \cref{rect-band}, the unique P-YBE solutions on $S$ are given by
	\begin{align*}
	r\left(a,b\right)=\left(\bar{\theta}\left(b\right), ab\right),
	\end{align*} 
where $\bar{\theta}$ is an idempotent endomorphism of $S$. \\
Since $S \in \mathcal{S}$ and $S$ is idempotent, by \cref{th-power}, we obtain that $r^4=r^2$ and the solutions $r^2$ and $r^3$ are respectively
	\begin{itemize}
	\item[-]$r^2\left(a,b\right)=\left( \bar{\theta}\left(a\right), \bar{\theta}\left(b\right)ab\right)$,
	\item[-]$r^3\left(a,b\right)=\left( \bar{\theta}\left(b\right), \bar{\theta}\left(a\right)ab\right)$.
	\end{itemize}
\end{ex}

\section{YBE solutions derived from PE solutions }

The aim of this section is to introduce a new method to construct solutions of the Yang-Baxter equation on the cartesian product of two sets through solutions of the pentagon equation.

\medskip

For the ease of the reader, before proving the main theorem we provide the conditions to obtain a YBE solution in an easier case involving a PE solution and an R-QYBE solution. 

Let $S$, $T$ be semigroups, $s$ a PE solution on $S$ given by $s(a,b)=(ab,\theta_a(b))$, and $t$ an R-QYBE solution on $T$ given by $t(u,v)=(\theta_v(u),vu)$. If $\alpha:T\to S^{S}$ is a map, set $\alpha_{u}:= \alpha\left(u\right)$, for every $u\in T$, and
\begin{align*}
&a_ub_v:=\alpha_u\left(a\right)\alpha_{\theta_v\left(u\right)}\left(b\right),
\end{align*}
for all $a,b \in S$ and $u,v \in T$.
If the following conditions hold
{\small
	\begin{center}
		\begin{minipage}[b]{.5\textwidth}
			\vspace{-\baselineskip}
			\begin{align}
			a\,b_uc_v&=a\theta_b\alpha_v(c)\,b_uc_v\label{p1m'}\tag{p1}
			\end{align}
		\end{minipage}%
	\end{center}
	\begin{center}
		\begin{minipage}[b]{.5\textwidth}
			\vspace{-\baselineskip}
			\begin{align}
			\theta_a\theta_b\alpha_u = \theta_{\alpha_v\left(b\right)}\alpha_{\theta_u\left(v\right)}\label{p2m'}\tag{p2}
			\end{align}
		\end{minipage}%
	\end{center}
	\begin{center}
		\begin{minipage}[b]{.5\textwidth}
			\vspace{-\baselineskip}
			\begin{align}\label{p3m'}\tag{p3}
			\thetaa{a}{(bc)}=\thetaa{a\theta_b\alpha_u(c)}{(bc)}
			\end{align}
		\end{minipage}%
	\end{center}
	\begin{center}
		\begin{minipage}[b]{.5\textwidth}
			\vspace{-\baselineskip}
			\begin{align}\label{p4m'}\tag{p4}
			a_ub_v=\alpha_{\theta_{wv}\left(u\right)}\left(a\alpha_v\left(b\right)\right)
			\end{align}
		\end{minipage}%
	\end{center}
	\begin{center}
		\begin{minipage}[b]{.5\textwidth}
			\vspace{-\baselineskip}
			\begin{align}\label{p5m'}\tag{p5}
			\thetaa{a}=\alpha_u\theta_a
			\end{align}
		\end{minipage}%
	\end{center}
}
\noindent for all $a,b,c \in S$ and $u,v,w \in T$, then the map
	\begin{align*}
r\left(a,u\, ; \,b,v\right)=\left(  \theta_a\alpha_u\left(b\right), vu \, ;\, a\alpha_u\left(b\right), \theta_v\left(u\right)\right),
\end{align*}
is a YBE solution on $S \times T$.\\
The proof of this statement is omitted since it will be contained in \cref{th:PE-YBE-1}. Let us see the efficacy of this method by using easy PE solutions.

\begin{ex}\label{ex-1}
Let $S \in \mathcal{W}=[abc=abdbc, \, a^3=a^2]$ (see \cite[p. 370]{Mo03}), $k$ a given element of $S$, and $s$ the PE solution on $S$ defined by $s(a,b)=(ab,k^2)$. Note that $s$ does not satisfy \eqref{YB1} and so it is not a QYBE solution. \\
	Moreover, let $T\in \mathcal{S}$ in \eqref{GRB}, and $t$ the R-QYBE solution on $T$ given by $t(u,v)=(u,vu)$. If we consider $\alpha_u(a)=k^2$, for every $a \in S$ and $u \in T$, then conditions from \eqref{p1m'} to \eqref{p5m'} are trivially satisfied, hence the map given by
	\begin{align*}
r\left(a,u\, ; \,b,v\right)=\left( k^2, vu\, ; \, ak^2, u\right)
\end{align*}
is a YBE solution on $S \times T$. One can check that $r^5=r^3$ and $r^2,r^3,r^4$ are still YBE solutions and they respectively are
\begin{itemize}
\item[-] $r^2\left(a,u\, ; \,b,v\right)=\left( k^2, uvu\, ; \,k^2, vu\right)$,
\item[-] $r^3\left(a,u\, ; \,b,v\right)=\left( k^2, vvu\, ; \,k^2, uvu\right)$,
\item[-] $r^4\left(a,u\, ; \,b,v\right)=\left( k^2, uvu\, ; \,k^2, vvu\right)$.
\end{itemize}
\end{ex}
\medskip

Note that if $S$ is a right zero semigroup, then by \eqref{omotheta} and \eqref{idemtheta} one can see that the unique PE solutions $s$ on $S$ are of the form 
$s\left(a,b\right) = \left(b, \varphi\left(b\right)\right)$, where $\varphi$ is an idempotent map from $S$ into itself. Clearly, such a map $s$ is a P-QYBE solution.
Thus, we provide the following example.
\begin{ex} \label{ex-2}
Let $S$ be a right zero semigroup, $\varphi$ an idempotent map from $S$ into itself, and $s$ the P-QYBE solution on $S$ given by $s\left(a,b\right) = \left(b, \varphi\left(b\right)\right)$. 
Moreover, let $T$ be a rectangular band and $t$ the R-QYBE solution on $T$ defined by $t(u,v)=(u,vu)$.\\ 
Set $\alpha_u:= \varphi$ from $S$ into itself, for every $u \in T$. Then, conditions from \eqref{p1m'} to \eqref{p5m'} trivially hold and so the map  
		\begin{align*}
		r\left(a,u\, ; \,b,v\right)
		= \left(\varphi(b), vu\, ; \,\varphi(b),u\right)
		\end{align*}
		is a YBE solution. Moreover, it holds $r^3=r$ and $r^2$ is still a YBE solution and it is given by
		\begin{align*}
		r^2\left(a,u\, ; \,b,v\right)
		= \left(\varphi(b), u\, ; \,\varphi(b),vu\right),
		\end{align*}
		for all $a,b \in S$ and $u,v \in T$.
\end{ex}
\medskip

In order to present the main theorem, we introduce the following definition. 
\begin{defin}
	Let $S$, $T$ be two semigroups, $s(a,b)=(ab,b)$ a PE solution on $S$ given by $s(a,b)=(ab,\theta_a(b))$, and $t$ a reversed solution on $T$ given by $t(u,v)=(\theta_v(u),vu)$. If $\alpha:T\to S^{S}$ and $\beta:S\to T^{T}$ are maps, set $\alpha_{u}:= \alpha\left(u\right)$, for every $u\in T$, and $\beta_{a}:= \beta\left(a\right)$, for every $a\in S$ and
	\begin{align*}
	&a_ub_v:=\alpha_u\left(a\right)\alpha_{\theta_v\beta_a\left(u\right)}\left(b\right) \qquad \theta_{ab_u}\left(c_v\right):=\theta_{a\alpha_u\left( b\right)}\alpha_{\theta_v\beta_b\left(u\right)}\left(c\right)\\
	&u_av_b:=\beta_a\left(u\right)\beta_{\theta_b\alpha_u\left(a\right)}\left(v\right)
	\qquad\theta_{uv_a}\left(w_b\right):=\theta_{u\beta_a\left(v\right)}\beta_{\theta_b\alpha_v\left(a\right)}\left(w\right),
	\end{align*}
	for all $a,b,c \in S$ and $u,v,w \in T$.
	If the following conditions hold
	{\small
		\begin{center}
			\begin{minipage}[b]{.5\textwidth}
				\vspace{-\baselineskip}
				\begin{align}
				a\,b_uc_v&=a\theta_b\alpha_v(c)\,b_uc_v\label{p1m}\tag{p1}
				\end{align}
			\end{minipage}%
			\hfill\hfill
			\begin{minipage}[b]{.5\textwidth}
				\vspace{-\baselineskip}
				\begin{align}\label{r1m}\tag{r1}
				u\, v_aw_b=u\theta_v\beta_b(w)\, v_aw_b
				\end{align}
			\end{minipage}
		\end{center}
		\begin{center}
			\begin{minipage}[b]{.5\textwidth}
				\vspace{-\baselineskip}
				\begin{align}
				\theta_a\theta_b\alpha_u = \theta_{\alpha_v\left(b\right)}\alpha_{\theta_u\beta_b\left(v\right)}\label{p2m}\tag{p2}
				\end{align}
			\end{minipage}%
			\hfill\hfill\hfill
			\begin{minipage}[b]{.5\textwidth}
				\vspace{-\baselineskip}
				\begin{align}\label{r2m}\tag{r2}
				\theta_u\theta_v\beta_a = \theta_{\beta_b\left(v\right)}\beta_{\theta_a\alpha_v\left(b\right)}
				\end{align}
			\end{minipage}
		\end{center}
		\begin{center}
			\begin{minipage}[b]{.5\textwidth}
				\vspace{-\baselineskip}
				\begin{align}\label{p3m}\tag{p3}
				\thetaa{a}{(bc)}=\thetaa{a\theta_b\alpha_u(c)}{(bc)}
				\end{align}
			\end{minipage}%
			\hfill\hfill
			\begin{minipage}[b]{.5\textwidth}
				\vspace{-\baselineskip}
				\begin{align}\label{r3m}\tag{r3}
				\thetaa{u}{(vw)}=\thetaa{u\theta_v\beta_a(w)}{(vw)}
				\end{align}
			\end{minipage}
		\end{center}
		\begin{center}
			\begin{minipage}[b]{.5\textwidth}
				\vspace{-\baselineskip}
				\begin{align}\label{p4m}\tag{p4}
				a_ub_v=\alpha_{\theta_{wv_b}\left(u_a\right)}\left(a\alpha_v\left(b\right)\right)
				\end{align}
			\end{minipage}%
			\hfill\hfill
			\begin{minipage}[b]{.5\textwidth}
				\vspace{-\baselineskip}
				\begin{align}\label{r4m}\tag{r4}
				u_av_b=\beta_{\theta_{cb_v}\left(a_u\right)}\left(u\beta_b\left(v\right)\right)
				\end{align}
			\end{minipage}
		\end{center}
		\begin{center}
			\begin{minipage}[b]{.5\textwidth}
				\vspace{-\baselineskip}
				\begin{align}\label{p5m}\tag{p5}
				\thetaa{a}=\alpha_u\theta_a
				\end{align}
			\end{minipage}%
			\hfill\hfill
			\begin{minipage}[b]{.5\textwidth}
				\vspace{-\baselineskip}
				\begin{align}\label{r5m}\tag{r5}
				\theta_{u}=\beta_a\theta_{u}
				\end{align}
			\end{minipage}
		\end{center}
	}
	\noindent for all $a,b,c \in S$ and $u,v,w \in T$, then $(s,t,\alpha,\beta)$ is called a \emph{pentagon quadruple}.
\end{defin}

\begin{theor}\label{th:PE-YBE-1}
	Let $S$, $T$ be semigroups and $(s,t,\alpha,\beta)$ a pentagon quadruple. Then, the map  $r: (S \times T)^2\to(S \times T)^2$ defined by
	\begin{align*}
	r\left(a,u\, ; \,b,v\right)=\left(\theta_a\alpha_u\left(b\right), v\beta_b\left(u\right) \, ; \, a\alpha_u\left(b\right), \theta_v\beta_b\left(u\right)\right),
	\end{align*}
	for	all $a,b \in S$ and $u,v \in T$, is a YBE solution on $S \times T$.
	\begin{proof}
		Let $(a,u), (b,v), (c,w) \in S \times T$. 	In order to show that $r$ is a YBE solution, set
		\begin{align*}
		\lambda_{\left(a,u\right)}\left(b,v\right):=\left(\theta_a\alpha_u\left(b\right), v\beta_b\left(u\right)\right)\qquad \rho_{\left(b,v\right)}\left(a,u\right):=\left(a\alpha_u\left(b\right), \theta_v\beta_b\left(u\right)\right),
		\end{align*}
	we have to check that
		\begin{align*}
		&\lambda_{\left(a,u\right)}\lambda_{\left(b,v\right)}\left(c,w\right) = \lambda_{\lambda_{\left(a,u\right)}\left(b,v\right)}\lambda_{\rho_{\left(b,v\right)}\left(a,u\right)}\left(c,w\right)\\
		&\rho_{\left(c,w\right)}\rho_{\left(b,v\right)}\left(a,u\right)
		= \rho_{\rho_{\left(c,w\right)}\left(b,v\right)}\rho_{\lambda_{\left(b,v\right)}\left(c,w\right)}\left(a,u\right)\\
		&\lambda_{\rho_{\lambda_{\left(b,vu\right)}\left(c,w\right)}\left(a,u\right)}\rho_{\left(c,w\right)}\left(b,v\right)
		= \rho_{\lambda_{\rho_{\left(b,v\right)}\left(a,u\right)}\left(c,w\right)}\lambda_{\left(a,u\right)}\left(b,v\right).
		\end{align*}
Initially, note that
\begin{align*}
 \lambda_{\left(a,u\right)}\lambda_{\left(b,v\right)}\left(c,w\right)=    \left(\theta_a\alpha_u\theta_b\alpha_v\left(c\right), w\beta_c\left(v\right)\beta_{\theta_b\alpha_v\left(c\right)}\left(u\right)\right),
\end{align*}
\begin{align*}
    &\lambda_{\lambda_{\left(a,u\right)}\left(b,v\right)}\lambda_{\rho_{\left(b,v\right)}\left(a,u\right)}\left(c,w\right)\\
    &=\left(\theta_{\theta_a\alpha_u\left(b\right)}\alpha_{v\beta_b\left(u\right)}\theta_{a\alpha_u\left(b\right)}\alpha_{\theta_v\beta_b\left(u\right)}\left(c\right), w \beta_c\theta_v\beta_b\left(u\right)\beta_{\theta_{a\alpha_u\left( b\right)}\alpha_{\theta_v\beta_b\left(u\right)}\left(c\right)}\left(v\beta_b\left(u\right)\right)\right),
\end{align*}
moreover
\begin{align*}
    &\rho_{\rho_{\left(c,w\right)}\left(b,v\right)}\rho_{\lambda_{\left(b,v\right)}\left(c,w\right)}\left(a,u\right)\\
    &=\left(a\alpha_u\theta_b\alpha_v\left(c\right)\alpha_{\theta_{w\beta_c\left( v\right)}\beta_{\theta_b\alpha_v\left(c\right)}\left(u\right)}\left(b\alpha_v\left(c\right)\right),
		\theta_{\theta_w\beta_c\left(v\right)}\beta_{b\alpha_v\left(c\right)}\theta_{w\beta_c\left(v\right)}\beta_{\theta_b\alpha_v\left(c\right)}\left(u\right)\right),
\end{align*}
\begin{align*}
 \rho_{\left(c,w\right)}\rho_{\left(b,v\right)}\left(a,u\right)=\left(a\alpha_u\left(b\right)\alpha_{\theta_v\beta_b\left(u\right)}\left(c\right),\theta_w\beta_c\theta_v\beta_b\left(u\right)\right), 
\end{align*}
and
\begin{align*}
    &\rho_{\lambda_{\rho_{\left(b,v\right)}\left(a,u\right)}\left(c,w\right)}\lambda_{\left(a,u\right)}\left(b,v\right)\\
    &=\left(\theta_a\alpha_u\left(b\right)\alpha_{v\beta_b\left(u\right)}\theta_{a\alpha_u\left(b\right)}\alpha_{\theta_v\beta_b\left(u\right)}\left(c\right), 
		\theta_{w \beta_c\theta_v\beta_b\left(u\right)} \beta_{\theta_{a\alpha_u\left(b\right)}\alpha_{\theta_v\beta_b\left(u\right)}\left(c\right)}   \left( v\beta_b\left(u\right)\right)\right),
\end{align*}
\begin{align*}
    &\lambda_{\rho_{\lambda_{\left(b,vu\right)}\left(c,w\right)}\left(a,u\right)}\rho_{\left(c,w\right)}\left(b,v\right)\\&=
    \left(\theta_{a\alpha_u\theta_b\alpha_v\left(c\right)}\alpha_{\theta_{w\beta_c\left( v\right)}\beta_{\theta_b\alpha_v\left(c\right)}\left(u\right)}\left(b\alpha_v\left(c\right)\right), \theta_w\beta_c\left(v\right)  \beta_{b\alpha_v\left(c\right)}\theta_{w\beta_c\left(v\right)}\beta_{\theta_b\alpha_v\left(c\right)}\left(u\right)\right).
\end{align*}
		For simplicity, we compare only the first components of each pair since the other equalities can be obtained reversing the role of the maps $\alpha$ and $\beta$. We have that 
\begin{align*}
    &\theta_a\alpha_u\theta_b\alpha_v\left(c\right) \\
		&=\theta_a\theta_b\alpha_v\left(c\right)&\mbox{by $\eqref{p5m}$}\\
		&=\theta_{\alpha_u\left(b\right)}\alpha_{\theta_v\beta_b\left(u\right)}\left(c\right)&\mbox{by $\eqref{p2m}$}\\
		&=\theta_{\theta_a\alpha_u\left(b\right)}\theta_{a\alpha_u\left(b\right)}\alpha_{\theta_v\beta_b\left(u\right)}\left(c\right)&\mbox{by \, $\eqref{idemtheta}$}\\
		&=\theta_{\theta_a\alpha_u\left(b\right)}\alpha_{v\beta_b\left(u\right)}\theta_{a\alpha_u\left(b\right)}\alpha_{\theta_v\beta_b\left(u\right)}\left(c\right),&\mbox{by  $\eqref{p5m}$}
\end{align*}
in addition
		\begin{align*} &a\alpha_u\theta_b\alpha_v\left(c\right)\alpha_{\theta_{w\beta_c\left( v\right)}\beta_{\theta_b\alpha_v\left(c\right)}\left(u\right)}\left(b\alpha_v\left(c\right)\right) \\
		&=a\theta_b\alpha_v\left(c\right)\alpha_{\theta_{wv_c}\left(u_b\right)}\left(b\alpha_v\left(c\right)\right)&\mbox{by $\eqref{p5m}$}\\
		&=a\theta_b\alpha_v\left(c\right)b_uc_v&\mbox{by $\eqref{p4m}$}\\\
		&=ab_uc_v&\mbox{by $\eqref{p1m}$}\\
		&=a\alpha_u\left(b\right)\alpha_{\theta_v\beta_b\left(u\right)}\left(c\right),
		\end{align*}
		and finally
		\begin{align*}
		&\theta_{a\alpha_u\theta_b\alpha_v\left(c\right)}\alpha_{\theta_{w\beta_c\left( v\right)}\beta_{\theta_b\alpha_v\left(c\right)}\left(u\right)}\left(b\alpha_v\left(c\right)\right)\\
		&=\theta_{a\theta_b\alpha_v\left(c\right)}\alpha_{\theta_{wv_c}\left(u_b\right)}\left(b\alpha_v\left(c\right)\right)&\mbox{by $\eqref{p5m}$}\\
		&=\theta_{a\theta_b\alpha_v\left(c\right)}\left(b_uc_v\right)&\mbox{by $\eqref{p4m}$}\\
		&=\theta_a\left(b_uc_v\right)&\mbox{by $\eqref{p3m}$}\\
		&=\theta_a\left(\alpha_u\left(b\right)\alpha_{\theta_v\beta_b\left(u\right)}\left(c\right)\right)\\
		&=\theta_a\alpha_u\left(b \right)\theta_{a\alpha_u\left(b\right)}\left(\alpha_{\theta_v\beta_b\left(u\right)}\left(c\right)\right)&\mbox{by \eqref{omotheta}}\\
		&=\theta_a\alpha_u\left(b\right)\alpha_{v\beta_b\left(u\right)}\theta_{a\alpha_u\left(b\right)}\alpha_{\theta_v\beta_b\left(u\right)}\left(c\right)&\mbox{by $\eqref{p5m}$}	\end{align*}
Therefore, $r$ is a YBE solution.
	\end{proof}
\end{theor}

\begin{ex}\label{ex-3}
	Let $S$ be a rectangular band and $s$ the P-QYBE solution on $S$ given by $s(a,b)=(ab,\gamma(b))$. Moreover, let $t$ be the Militaru solution on a semigroup $T$ given by $t(u,v)=(f(u),f(v))$.\\
	Set $\alpha_u=\gamma$, for every $u \in T$, and $\beta_a=f$, for every $a \in S$. Then, $(s,t,\alpha,\beta)$ is a pentagon quadruple and, by \cref{th:PE-YBE-1}, the map
	\begin{align*}
	r\left(a,u\, ; \, b,v\right)=
	\left(\gamma\left(b\right),  f(v)\, ; \, a\gamma(b), f\left(u\right) \right)
	\end{align*}
	is a YBE solution on $S \times T$. Moreover, it holds $r^4=r^2$ and the powers of $r$ are still YBE solutions and they are 
	\begin{itemize}
	    \item[-] $r^2\left(a,u\, ; \, b,v\right)=
	\left(\gamma\left(ab\right),  f(u)\, ; \,  \gamma(b), f\left(v\right)\right)$
	    \item[-] $r^3\left(a,u\, ; \, b,v\right)=
	\left(\gamma\left(b\right),  f(v)\, ; \,  \gamma(ab), f\left(u\right)\right)$.
	\end{itemize}
\end{ex}
\medskip

\begin{ex}
Assume that $S$ and $T$ are groups and $(s,t,\alpha, \beta)$ is a pentagon quadruple. Conditions \eqref{p1m} and \eqref{r1m} become 
\begin{align}
	\theta_a\alpha_u(b)&=1_S\label{marzia}\\ 
	\theta_u\beta_a(v)&=1_T\label{paola}.
\end{align}
By \eqref{marzia}, \eqref{paola}, and \cite[Lemma 11]{CaMaMi19x}, it follows that \eqref{p2m}, \eqref{r2m}, \eqref{p3m}, and \eqref{r3m} are trivially  satisfied.
Moreover, \eqref{p4m} and \eqref{r4m} become
\begin{align}
			\alpha_u(a)\alpha_{1_T}(b)&=\alpha_{1_T}(a\alpha_v(b)) \label{8}\\ \beta_a(u)\beta_{1_S}(v)&=\beta_{1_S}(u\beta_b(v))\label{9}
\end{align}
By \eqref{p5m}, we get $\alpha_u(1_S)=1_S$, for every $u \in T$ and similarly, by \eqref{r5m}, $\beta_a(1_T)=1_T$, for every $a \in S$. Moreover, by \eqref{8}, with $b=1_S$ and $v=1_T$ we obtain 
\begin{center}
 $\alpha_u(a)=\alpha_{1_T}(a)$,    
\end{center}
for all $u \in T$ and $a \in S$. Thus, all the maps $\alpha_u$ are all equal to a map $\bar{\alpha}$, for every $u \in T$. Similarly, the maps $\beta_a$ are all equal to $\bar{\beta}$, for all $a \in S$. Therefore, by \cref{th:PE-YBE-1}, the map
\begin{align*}
	r\left(a,u\, ; \,b,v\right)=\left(1_S, v\bar{\beta}\left(u\right) \, ; \, a\bar{\alpha}\left(b\right), 1_T\right),
	\end{align*}
is a YBE solution on $S \times T$. Note that 
\begin{align*}
	r^n\left(a,u\, ; \,b,v\right)=\left(1_S, \bar{\beta}^{n-1}\left(v\right)\bar{\beta}\left(u\right) \, ; \, \bar{\alpha}^{n-1}\left(a\right)\bar{\alpha}\left(b\right), 1_T\right),
	\end{align*}
for every $n\in \mathbb{N}$. In this way, by choosing $\bar{\alpha}$ an idempotent endomorphism of $S$ and $\bar{\beta}$ an idempotent endomorphism of $T$, then we obtain examples of solutions $r$ with the property $r^3=r^2$.\\
To find other examples on groups, by conditions \eqref{marzia} and \eqref{paola} one has to look at maps $\bar{\alpha}$ and $\bar{\beta}$ whose images $\im\bar{\alpha}$ and $\im\bar{\beta}$ are contained in the kernel $K_{s}$ of $s$ and $K_{t}$ of $t$, respectively. Given a PE solution $s(a,b) = (ab, \theta_a(b))$ on a group $S$, we remind that the subset of $S$ defined by
\begin{align*}
    K_{s}:= \lbrace a \, | \, a\in S, \ \theta_{1_S}(a) = 1_S\rbrace
\end{align*}
is a normal subgroup of $S$ called the \emph{kernel of $s$} (cf. \cite[Lemma 13]{CaMaMi19x}).
\end{ex}
\bigskip

\section{Some comments and questions}

We conclude with some remarks and questions that arise throughout this work. Specifically, our purpose is to show that if some properties hold for a YBE solution $r$, then it does not necessarily hold for the QYBE solution $s:=\tau r$ and vice versa.
\medskip

First, we give the definitions of the index and the period of a map $f$ as 
\begin{align*}
        \indd{\left(f\right)}
        &:=\min\left\{\left.j \,\right|\, j\in\mathbb{N}_0, \, \exists \, l\in \mathbb{N}\ f^j = f^l , \ j\neq l\right\}\\
        \perr{\left(f\right)}
        &:=\min\left\{\left.k\,\right| \, k\in\mathbb{N}, \, f^{\indd{\left(r\right)}+k} = f^{\indd{\left(r\right)}}\right\},
    \end{align*}
respectively. As observed in \cite{CaCoSt19}, these definitions are slightly different from the classical ones (cf. \cite[p. 10]{Ho95}).
\medskip

\noindent In general, we note that the index or the period of a QYBE solution $s$ are independent from the index or the period of its YBE solution $r$. The following are examples of QYBE solutions $s$ for which it holds $\indd(s)=1$ and $\perr(s)=1$, but the corresponding YBE solutions $r=\tau s$ have different periods.

\begin{exs}\label{exsopra}\hspace{1mm}
\begin{enumerate}
    \item Let $S$ be a rectangular band, i.e., an idempotent semigroup with the property $abc=ac$, for all $a,b,c \in S$, and $s$ the P-QYBE solution on $S$ given by $s(a,b)=(ab,b)$. Then, they hold $\indd{(r)}=1$ and $\perr{(r)}=2$.
    
    \vspace{1mm}
    
    \item Let $S$ be a monoid and $s$ the solution P-QYBE on $S$ given by $s(a,b)=(ab,1)$. Then, in this case we have $\indd{(r)}=1$ and $\perr{(r)}=1$.
    \end{enumerate}
\end{exs}
\medskip

\noindent Another aspect is that if $r$ is a YBE solution for which its powers $r^n$ are still solutions, then the maps $s^n:=(\tau r)^n$ are not necessarily QYBE solutions. The following is an example of YBE solution $r$ for which the powers are all solutions while the powers of $s= \tau r$ are not all QYBE solutions.
\begin{ex} 
Let $r$ be the P-YBE solution given by $r(a,b) = (\gamma(b), ab)$ on a semigroup $S\in \mathcal{S}$. Then, by \cref{th-power}, we have that $r^5=r^3$ and the powers of $r$ are still YBE solutions. 
Now, if we consider $s=\tau r$, it holds that $s^4 = s^3$ and one can check that $s^3$ is still a QYBE solution, while $s^2$ is not. 
\end{ex}

\noindent As observed in \cref{ex:powers}.2, the powers of the YBE solution $r(a,b)=(bab, ab)$ are all YBE solutions and as seen in \cref{ex:powers-2}.1, the map $s=\tau r$ is such that its powers are all Q-YBE solutions.
Hence, a question arises.
\vspace{3mm}

\noindent \textbf{Question.} Classify the YBE solutions $r$ such that their powers are still YBE solutions and for which also the powers of the QYBE solution $s= \tau r$ are still QYBE solutions.

\bigskip


\bibliography{bibliography}

\end{document}